\documentclass[]{interact}

\usepackage{epstopdf}
\usepackage{graphicx}
\usepackage{booktabs}
\usepackage{placeins}
\usepackage{enumitem}
\usepackage{url}
\usepackage[numbers,sort&compress]{natbib}
\bibpunct[, ]{[}{]}{,}{n}{,}{,}
\renewcommand\bibfont{\fontsize{10}{12}\selectfont}

\theoremstyle{plain}
\newtheorem{theorem}{Theorem}[section]
\newtheorem{lemma}[theorem]{Lemma}
\newtheorem{corollary}[theorem]{Corollary}

\theoremstyle{definition}
\newtheorem{definition}[theorem]{Definition}
\newtheorem{example}[theorem]{Example}
\theoremstyle{remark}
\newtheorem{remark}[theorem]{Remark}

\newcommand{\R}{\mathbb{R}}
\newcommand{\xb}{\bar{x}}
\newcommand{\Ec}{E_c}
\newcommand{\Ecb}{\overline{E}_c}
\newcommand{\T}{\mathsf{T}}

\begin{document}

\articletype{ARTICLE}

\title{Stability of equilibria in discrete-time modified-gradient systems and the multivariate breeder's equation}

\author{
\name{Benjamin J. Ridenhour\textsuperscript{a,b,c}\thanks{CONTACT B.~J. Ridenhour. Email: bridenhour@uidaho.edu}}
\affil{\textsuperscript{a}Department of Mathematics and Statistical Science, University of Idaho, Moscow, ID, USA; \textsuperscript{b}Bioinformatics and Computational Biology Program, University of Idaho, Moscow, ID, USA; \textsuperscript{c}Institute for Interdisciplinary Data Science, University of Idaho, Moscow, ID, USA}
}

\maketitle

\begin{abstract}
The multivariate breeder's equation $\Delta\bar z=G\nabla\ln\bar W$ is a widely used discrete-time modified-gradient system describing the evolution of continuous traits in organisms.  This system could more conventionally be expressed as $x_{k+1}=x_k+P_k\nabla f(x_k)$. For its continuous-time counterpart $x'=P(t)\nabla f(x)$, it has been proven that an isolated local maximum of $f$ is uniformly stable whenever $P(t)$ is positive semidefinite, as is the case for the ``$G$-matrix.'' Furthermore, the continuous time system is uniformly asymptotically stable when the maximum is an isolated critical point and $\int_0^\infty\lambda_1(P(t))\,dt=\infty$. We show that the literal discrete analog of this result is false: with constant positive definite $P$ the iterates can oscillate or diverge. The result is restored by a single additional hypothesis, a uniform bound on the largest eigenvalue of $P_k$, $\lambda_n$, such that $\lambda_n(P_k)\le 2/L$  where $L$ is a local Lipschitz constant of $\nabla f$. Under this bound the maximum is uniformly stable; under a stricter version together with $\sum_k\lambda_1(P_k)=\infty$, the maximum is uniformly asymptotically stable. Given these hypotheses, a level-set estimate of the basin of attraction carries over. Both bounds are sharp, and each remaining hypothesis is shown to be necessary. Biologically, the bound says that genetic variance in every direction must be small relative to the curvature of log mean fitness; outside this regime the mean phenotype can overshoot a fitness peak as hypothesized by Fisher's geometric model. As further application of our results, the same conditions apply to the mean update of Gaussian evolutionary optimization algorithms, in which the covariance of the search population plays the role of $G$.
\end{abstract}


\begin{amscode}
39A30; 92D15; 37C75; 90C30
\end{amscode}

\section{Introduction}

Quantitative genetic models of multivariate phenotypic evolution describe the change in a vector of mean continuous traits $\bar z$ by Lande's equation
\begin{equation}\label{eq:lande}
  \Delta\bar z \;=\; G\,\nabla\ln\bar W(\bar z),
\end{equation}
where $G$ is the additive genetic variance--covariance matrix and $\bar W$ is mean fitness~\cite{Lush1937,Lande1976,Lande1979,LandeArnold1983}. The $G$-matrix is symmetric positive semidefinite, and it is now well documented that it too changes over evolutionary time~\cite{Arnold2008}. Equation~\eqref{eq:lande} is therefore an instance of a \emph{modified-gradient system}: the trait vector climbs the surface described by $\ln\bar W$, but in a direction deflected by a time-varying positive semidefinite matrix.

In earlier work~\cite{RR2018}, we studied the continuous-time version
\begin{equation}\label{eq:E}
  x' \;=\; P(t)\,\nabla f(x),
\end{equation}
with $P(t)$ symmetric positive semidefinite, and asked when a local maximum $\xb$ of $f$ is a stable equilibrium regardless of how $P(t)$ evolves. The answer was that an isolated local maximum is always uniformly stable, and is uniformly asymptotically stable provided it is also an isolated critical point and the \emph{eigenvalue condition}
\begin{equation}\label{eq:EC}
  \int_0^\infty\lambda_1\bigl(P(t)\bigr)\,dt=\infty
\end{equation}
holds, where $\lambda_1$ denotes the smallest eigenvalue of $P(t)$. A level-set estimate of the basin of attraction followed under further conditions on $f$. Biologically, these results say that the evolution of $G$ cannot destabilize a fitness peak; rather, changes in the genetic covariance matrix can at most it can slow convergence to it in continuous time.

Many organisms, however, have discrete, non-overlapping generations, and~\eqref{eq:lande} is itself a per-generation recursion. The natural object of study for such evolutionary systems is then the difference equation
\begin{equation}\tag{$\Delta$}\label{eq:Delta}
  x_{k+1} \;=\; x_k + P_k\,\nabla f(x_k), \qquad k=0,1,2,\dots,
\end{equation}
with $P_k=G_k$ and $f=\ln\bar W$ from~\eqref{eq:lande}. The question addressed here is whether the conclusions of~\cite{RR2018} remain true for~\eqref{eq:Delta}.

They do not, at least not as stated. In continuous time, a trajectory cannot cross a level set of $f$ without passing through it, and it cannot pass through the equilibrium $\xb$ at all \cite{RR2018}. However, the discrete time step $P_k\nabla f(x_k)$ has finite length and can overshoot the optimum. A one-dimensional quadratic example (Section~\ref{sec:fail}) shows that, with every hypothesis of~\cite{RR2018} enforced, the iterates of~\eqref{eq:Delta} may converge, oscillate, or diverge, depending only on the product of the step size and the curvature of $f$.

Our main results show restricting the step size is the only obstruction. Specifically, if $\nabla f$ is $L$-Lipschitz near $\xb$, then a bound on the largest eigenvalue of $P_k$, $\lambda_n$, such that $\lambda_n(P_k)\le2/L$  restores uniform stability (Theorem~\ref{thm:stab}(i)); the strict bound $\lambda_n(P_k)\le(2-\eta)/L$, together with the discrete eigenvalue condition $\sum_k\lambda_1(P_k)=\infty$, restores uniform asymptotic stability (Theorem~\ref{thm:stab}(ii)); and the basin of attraction estimate of~\cite{RR2018} carries over under the same bounds (Theorem~\ref{thm:basin}). Both bounds are sharp, and each of the other underlying hypotheses are required (Section~\ref{sec:necessity}). The key concept is developing a one-step version of the trajectory derivative $V'=-\nabla f\cdot P\nabla f$ (Lemma~\ref{lem:ascent}) for a Lyapunov function $V(x)$, which is the matrix-step form of the classical descent lemma from optimization problems~\cite{Polyak1987,Bertsekas1999,Nesterov2004}.

The paper is organized as follows. Section~\ref{sec:setting} fixes notation and hypotheses. Section~\ref{sec:fail} gives the counterexample to the literal analog. Sections~\ref{sec:lemma}, ~\ref{sec:stab}, and ~\ref{sec:basin} prove the key lemma, the stability theorems, and the basin of attraction, respectively. Section~\ref{sec:necessity} shows that underlying hypotheses cannot be ignored. Lastly, Section~\ref{sec:discussion} discusses these results in the light of biological meaning, optimization, and some open problems.

\section{Setting and hypotheses}\label{sec:setting}

Notation follows~\cite{RR2018}. $|\cdot|$ is the Euclidean norm, $B_r(\xb)$ the open ball of radius $r$ about $\xb$, and $\overline{B}_r(\xb)$ its closure. For a symmetric matrix $S$, $\lambda_1(S)\le\cdots\le\lambda_n(S)$ are its eigenvalues. For symmetric $A,B$, $A\preceq B$ means that $B-A$ is positive semidefinite (the L\"owner order~\cite[\S7.7]{HornJohnson2013}). Throughout, $\xb$ is the equilibrium of interest, $M=f(\xb)$, and $V(x)=M-f(x)$.

We use the following hypotheses. H1 and H2 are those of~\cite{RR2018}; H0$_\Delta$ and H3$_\Delta$ are discrete versions of H0 and H3 of~\cite{RR2018}; the step-size conditions S1--S2 have no continuous-time counterpart.
\begin{description}[leftmargin=2.8em,style=nextline]
  \item[H0$_\Delta$] $D\subset\R^n$ is a domain, $f\in C^1(D)$, and $\nabla f$ is locally Lipschitz on $D$. Each $P_k$ is a symmetric positive semidefinite $n\times n$ matrix.
  \item[H1] $f$ has an isolated local maximum at $\xb\in D$.
  \item[H2] $\xb$ is an isolated critical point of $f$.
  \item[H3$_\Delta$] $\sum_{k=0}^\infty\lambda_1(P_k)=\infty$.
\end{description}
The step-size conditions refer to constants $r_0>0$, $L>0$ and $\Lambda>0$ such that $\overline{B}_{r_0}(\xb)\subset D$,
\[
  |\nabla f(u)-\nabla f(v)|\le L|u-v|\quad\text{for }u,v\in\overline{B}_{r_0}(\xb),
  \qquad\text{and}\qquad
  \lambda_n(P_k)\le\Lambda\quad\text{for all }k.
\]
\begin{description}[leftmargin=2.8em]
  \item[S1] $L\Lambda\le2$.
  \item[S2] $L\Lambda\le2-\eta$ for some $\eta\in(0,2]$.
\end{description}
H0 of~\cite{RR2018}, which asks that $\nabla f$ have $C^1$ components, implies the Lipschitz condition in H0$_\Delta$. 

Additional basin of attraction conditions (H4--H6; S3) are stated in Section~\ref{sec:basin}.

\begin{definition}
A sequence of~\eqref{eq:Delta} started at index $k_0$ \emph{exists} as long as its iterates remain in $D$. The equilibrium $\xb$ is \emph{stable} if for every $\varepsilon>0$ there are $\delta>0$ and $k_\varepsilon\ge0$ such that $|x_{k_0}-\xb|<\delta$ and $k_0\ge k_\varepsilon$ imply $x_k\in D$ and $|x_k-\xb|<\varepsilon$ for all $k\ge k_0$; the sequence is \emph{asymptotically stable} if in addition $x_k\to\xb$. Stability is \emph{uniform} if one can take $k_\varepsilon=0$, so that $\delta$ does not depend on the starting index.
\end{definition}
These are the discrete counterparts of the definitions of Hartman~\cite{Hartman1964} used in~\cite{RR2018}; for Lyapunov methods for difference equations see~\cite{KalmanBertram1960,LaSalle1986,Elaydi2005}.

\section{Failure of the discrete analog $\Delta$}\label{sec:fail}

\begin{example}\label{ex:quad}
Let $n=1$, $f(x)=-\tfrac a2x^2$ with $a>0$, and $P_k\equiv p>0$. Then $\xb=0$ is an isolated maximum and an isolated critical point. Furthermore, $\sum_k\lambda_1(P_k)= \sum_kp=\infty$, and $\nabla f$ is Lipschitz with $L=a$. By substitution, system~\eqref{eq:Delta} is then $x_{k+1}=(1-pa)x_k$; this recursion subsequently telescopes to $x_k=(1-pa)^{k-k_0}x_{k_0}$. From this, we can clearly see the behavior of the dynamical system depends only on $pa (=L\Lambda)$ (Table \ref{tab:quad}).

\begin{table}
\tbl{Dynamics of Example~\ref{ex:quad}.\label{tab:quad}}
{\begin{tabular}{@{}lll@{}}
\toprule
$pa=L\Lambda$ & Behavior of $x_k$ & Dynamics\\
\midrule
$0<pa<2$ & $x_k\to0$ geometrically & Uniformly asymptotically stable\\
$pa=2$ & $x_k=(-1)^{k-k_0}x_{k_0}$ & Uniformly stable, not asymptotically stable\\
$pa>2$ & $|x_k|\to\infty$ & Unstable\\
\bottomrule
\end{tabular}}
\end{table}
\end{example}

\begin{figure}
\centering
\includegraphics[width=\textwidth]{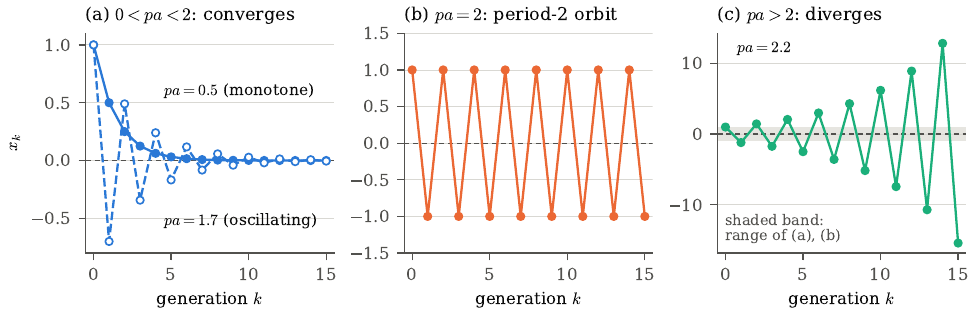}
\caption{Trajectories of Example~\ref{ex:quad}. For this example, $x_{k+1}=(1-pa)x_k$, where $\bar x = 0$ and $x_0 = 1$; the dashed horizontal line marks the maximum $\bar{x}=0$. (a) For $0<pa<2$ the iterates converge: monotonically when $pa<1$ (solid, $pa=0.5$) and with alternating overshoot when $1<pa<2$ (dashed, $pa=1.7$). (b) At $pa=2$ each step overshoots by exactly the starting distance, giving a period-2 orbit. (c) For $pa>2$ each iterate lands farther from the maximum than it started and the iterates diverge ($pa=2.2$); the shaded band $-1\le x_k\le1$ shows the vertical range of panels (a) and (b).}
\label{fig:quad}
\end{figure}

All hypotheses of Theorem~2.1 of~\cite{RR2018} hold in every row of Table \ref{tab:quad}, so neither part of that theorem transfers to~\eqref{eq:Delta} without an additional assumption. The thresholds $pa\le2$ and $pa<2$ are exactly S1 and S2, so the bounds in Theorem~\ref{thm:stab} cannot be improved. The mechanism is overshoot (Figure~\ref{fig:quad}): for $pa>1$ each step crosses the maximum, and for $pa>2$ it lands farther away than it started.

\section{Discrete Lyapunov decrease}\label{sec:lemma}

The following lemma replaces the identity $V'=-\nabla f\cdot P\nabla f$ that drives the continuous-time proofs. It is the descent lemma for functions with Lipschitz gradient~\cite{Polyak1987,Bertsekas1999},~\cite[Lemma~1.2.3]{Nesterov2004}, adapted to a matrix step.

\begin{lemma}[One-step ascent]\label{lem:ascent}
Let $P$ be symmetric positive semidefinite, $g=\nabla f(x)$ and $d=Pg$, and suppose $\nabla f$ is $L$-Lipschitz on a set containing the segment $[x,x+d]$. Then for every $s\in[0,1]$,
\begin{equation}\label{eq:ascent}
  f(x+sd)\;\ge\;f(x)+s\Bigl(1-\tfrac s2L\lambda_n(P)\Bigr)\,g^\T Pg.
\end{equation}
In particular, if $\lambda_n(P)\le\Lambda$ and $L\Lambda\le2$, then
\begin{equation}\label{eq:decrease}
  V(x+d)\;\le\;V(x)-\Bigl(1-\tfrac{L\Lambda}2\Bigr)g^\T Pg\;\le\;V(x)-\Bigl(1-\tfrac{L\Lambda}2\Bigr)\lambda_1(P)\,|g|^2.
\end{equation}
\end{lemma}

\begin{proof}
Let $\varphi(\tau)=f(x+\tau d)$ for $\tau\in[0,1]$. Then $\varphi\in C^1$ with $\varphi'(\tau)=\nabla f(x+\tau d)\cdot d$; in particular $\varphi'(0)=g\cdot d=g^\T d=g^\T Pg$. By Cauchy--Schwarz and the Lipschitz bound,
\[
  |\varphi'(\tau)-\varphi'(0)|\le|\nabla f(x+\tau d)-\nabla f(x)|\,|d|\le L\tau|d|^2 .
\]
Integrating from $0$ to $s$,
\[
  \varphi(s)\;\ge\;\varphi(0)+s\,\varphi'(0)-\tfrac{L s^2}{2}|d|^2
  \;=\;f(x)+s\,g^\T Pg-\tfrac{Ls^2}{2}\,g^\T P^2g .
\]
Every eigenvalue $\mu$ of $P$ satisfies $0\le\mu^2\le\lambda_n(P)\mu$, so $P^2\preceq\lambda_n(P)P$ and $g^\T P^2g\le\lambda_n(P)\,g^\T Pg$. This gives~\eqref{eq:ascent}. Taking $s=1$ and using $\lambda_n(P)\le\Lambda$ gives the first inequality in~\eqref{eq:decrease}; the second follows from $g^\T Pg\ge\lambda_1(P)|g|^2$ because the coefficient $1-L\Lambda/2$ is nonnegative.
\end{proof}

We record three consequences, each valid as long as the step $[x_k,x_{k+1}]$ stays where the Lipschitz bound holds.
\begin{enumerate}[label=(C\arabic*)]
  \item\label{c:mono} Under S1, $V(x_{k+1})\le V(x_k)$: the discrete analog of $V'\le0$.
  \item\label{c:strict} Under S2, $V(x_{k+1})\le V(x_k)-\tfrac\eta2\lambda_1(P_k)|\nabla f(x_k)|^2$: the analog of $V'\le-\lambda_1(P(t))|\nabla f|^2$.
  \item\label{c:segment} Under S1, $f(x_k+sd_k)\ge f(x_k)$ for every $s\in[0,1]$, because $s(1-sL\Lambda/2)\ge0$. Thus $f$ does not fall below its starting value anywhere along a step.
\end{enumerate}
We also use the step-length bound $|P\nabla f(x)|\le\lambda_n(P)|\nabla f(x)|\le\Lambda|\nabla f(x)|$, which keeps each segment inside the Lipschitz bound.

\section{Stability and asymptotic stability}\label{sec:stab}

\begin{theorem}\label{thm:stab}\leavevmode
\begin{enumerate}[label=(\roman*)]
  \item If H0$_\Delta$, H1 and S1 hold, then $\xb$ is a uniformly stable equilibrium of~\eqref{eq:Delta}.
  \item If H0$_\Delta$, H1, H2, H3$_\Delta$ and S2 hold, then $\xb$ is a uniformly asymptotically stable equilibrium of~\eqref{eq:Delta}.
\end{enumerate}
\end{theorem}

\begin{proof}[Proof of (i)]
Let $\varepsilon>0$. By H1 choose $r\in(0,\min(\varepsilon,r_0)]$ such that $f(x)<M$ for $0<|x-\xb|\le r$. In continuous time, a trajectory cannot reach the sphere $|x-\xb|=r$ without passing through every intermediate level of $V$.  In contrast, discrete time iterates can jump, so we need to build a barrier of positive thickness.

Let $\Gamma(\rho)=\max\{|\nabla f(x)|:|x-\xb|\le\rho\}$. Then $\Gamma$ is continuous and nondecreasing with $\Gamma(0)=0$, so we may choose $\rho\in(0,r)$ such that $\Lambda\Gamma(\rho)<r-\rho$ (i.e., the step-length bound $\Lambda |\nabla f(x)|$ prohibits crossing $r$). On the compact annulus $A=\{\rho\le|x-\xb|\le r\}$ the function $V$ is continuous and positive; let $m=\min_AV>0$. Finally choose $\delta\in(0,\rho]$ with $V<m$ on $B_\delta(\xb)$.

We claim that $|x_{k_0}-\xb|<\delta$ implies $x_k\in B_\rho(\xb)$ and $V(x_k)<m$ for all $k\ge k_0$. This holds for $k=k_0$. If it holds for $k$, then $|x_{k+1}-x_k|\le\Lambda\Gamma(\rho)<r-\rho$, so the segment $[x_k,x_{k+1}]$ lies in the convex set $B_r(\xb)$, where $\nabla f$ is $L$-Lipschitz. By~\ref{c:mono}, $V(x_{k+1})\le V(x_k)<m$. Because $V\ge m$ on $A$ and $x_{k+1}\in B_r(\xb)$, it follows that $x_{k+1}\in B_\rho(\xb)$.

Thus $|x_k-\xb|<\rho<\varepsilon$ for all $k\ge k_0$, and because $\delta$ does not depend on $k_0$ the stability is uniform.
\end{proof}

\begin{proof}[Proof of (ii)]
Repeat the construction for (i), however, given H2 we have $\nabla f(x)\ne0$ for $0<|x-\xb|\le r$ (i.e., $\bar x$ is an isolated critical point). Also let $|x_{k_0}-\xb|<\delta$. If $x_k=\xb$ for some $k$, the sequence remains at $\xb$; unlike in continuous time this can occur in finite time (see Figure \ref{fig:quad}a). Otherwise $x_k\in B_\rho(\xb)\setminus\{\xb\}$ for all $k\ge k_0$.

The sequence $V_k=V(x_k)$ is nonincreasing and nonnegative, so $V_k\to c\ge0$. Suppose $c>0$, and choose $\delta_1\in(0,\rho)$ with $V<c$ on $B_{\delta_1}(\xb)$. Then every iterate lies in the compact set $K=\{\delta_1\le|x-\xb|\le\rho\}$, on which $m_1=\min_K|\nabla f|^2>0$ by H2. By~\ref{c:strict},
\[
  V_k\;\le\;V_{k_0}-\frac{\eta m_1}{2}\sum_{j=k_0}^{k-1}\lambda_1(P_j)\;\longrightarrow\;-\infty,
\]
by H3$_\Delta$, a contradiction. Hence $V_k\to0$. If $x_k\not\to\xb$, there is $\varepsilon_1\in(0,\rho)$ with $\varepsilon_1\le|x_k-\xb|<\rho$ for infinitely many $k$; on this annulus $V$ has a positive minimum, contradicting $V_k\to0$.
\end{proof}

Here we show that the uniform margin $\eta$ in S2 can be replaced by a weighted divergence condition and asymptotic stability is still achieved.

\begin{corollary}\label{cor:weighted}
Suppose H0$_\Delta$, H1, H2 and S1 hold, and
\begin{equation}\label{eq:weighted}
  \sum_{k=0}^\infty\lambda_1(P_k)\bigl(2-L\lambda_n(P_k)\bigr)=\infty .
\end{equation}
Then $\xb$ is uniformly asymptotically stable.
\end{corollary}

\begin{proof}
Under S2 the weights satisfy $2-L\lambda_n(P_k)\ge\eta$, so~\eqref{eq:weighted} follows from H3$_\Delta$; the corollary therefore generalizes Theorem~\ref{thm:stab}(ii). The proof is that of Theorem~\ref{thm:stab}(ii), except that Lemma~\ref{lem:ascent} is applied with $\lambda_n(P_k)$ in place of $\Lambda$. With $g_k=\nabla f(x_k)$ this gives
\[
  V_{k+1}\;\le\;V_k-\tfrac12\,w_k\,|g_k|^2,\qquad w_k=\lambda_1(P_k)\bigl(2-L\lambda_n(P_k)\bigr)\ge0,
\]
and the rest of the argument goes through with $\eta\lambda_1(P_j)$ replaced by $w_j$.
\end{proof}

\begin{figure}
\centering
\includegraphics[width=0.7\textwidth]{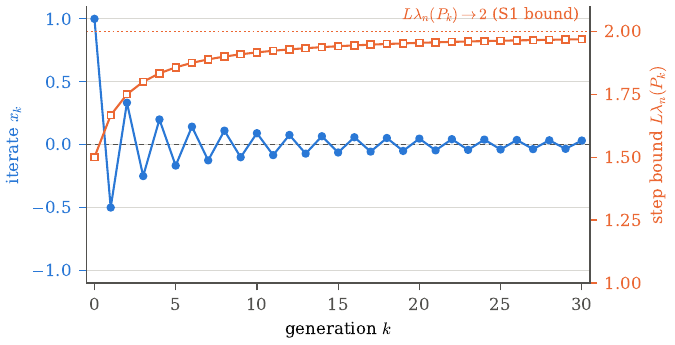}
\caption{Iterates of Example~\ref{ex:nomargin} with $x_0=1$. Left axis (circles): $x_k=(-1)^k/(k+1)$, which converges to the maximum $\bar x=0$ (dashed line) while alternating in sign. Right axis (squares): the step bound $L\lambda_n(P_k)=a\,p_k=2-(k+2)^{-1}$, which increases toward the S1 bound $L\lambda_n(P_k)=2$ (dotted line), so S2 fails for every $\eta>0$.}
\label{fig:nomargin}
\end{figure}

\begin{example}\label{ex:nomargin}
Take $f(x)=-\tfrac a2x^2$ and $p_k=\bigl(2-(k+2)^{-1}\bigr)/a$. Then $Lp_k\to2$, so S2 fails for every $\eta>0$. However, noting that $\tfrac 3{2a}\le p_k\le \tfrac 2a$ shows that $w_k=p_k/(k+2)\ge\tfrac32/(a(k+2))$ and~\eqref{eq:weighted} holds because the harmonic series diverges. Indeed $x_{k+1}=-\tfrac{k+1}{k+2}x_k$, so $|x_k|=|x_0|/(k+1)\to0$, with alternating sign (Figure~\ref{fig:nomargin}). This example shows that S2 is sufficient but not necessary for asymptotic stability.
\end{example}

\begin{remark}
H3$_\Delta$ plays the role of the classical step-size condition $\sum_k\alpha_k=\infty$ for gradient methods~\cite{Bertsekas1999,BertsekasTsitsiklis2000}, where it is paired with an upper bound such as $\alpha_k\le(2-\eta)/L$, the scalar form of S2. Theorem~\ref{thm:stab} is the matrix-step, local-stability counterpart, with $\lambda_1(P_k)$ and $\lambda_n(P_k)$ in the roles of the lower and upper step sizes.
\end{remark}

\section{Basin of attraction}\label{sec:basin}

As in Theorem~3.1 of~\cite{RR2018}, fix $c<M$, let $O_c=\{\xb\}\cup\{x:c<f(x)<M\}$, and let $\Ec$ be the component of $O_c$ containing $\xb$, with closure $\Ecb$ and boundary $\partial\Ec$. The hypotheses of~\cite{RR2018} are:
\begin{description}[leftmargin=2.8em]
  \item[H4] $\Ecb$ is bounded and $\Ecb\subset D$;
  \item[H5] $f(x)=c$ for all $x\in\partial\Ec$;
  \item[H6] $f$ has no critical point in $\Ecb$ other than $\xb$.
\end{description}
The discrete version of this proof needs one more hypothesis than the continuous version:
\begin{description}[leftmargin=2.8em]
  \item[S3] There are $L,\Lambda,\eta>0$ such that $\nabla f$ is $L$-Lipschitz on $\Ecb$, $\lambda_n(P_k)\le\Lambda$ for all $k$, and $L\Lambda\le2-\eta$.
\end{description}
By H4, $\Ecb$ is a compact subset of $D$, and a locally Lipschitz map is Lipschitz on compact sets, so by H0$_\Delta$ some $L$ always exists; S3 only asks that $\lambda_n(P_k)$ be uniformly below $2/L$.

\begin{theorem}\label{thm:basin}
Assume H0$_\Delta$, H1, H2, H3$_\Delta$, H4--H6 and S3. Then every sequence of~\eqref{eq:Delta} with $x_{k_0}\in\Ec$ exists for all $k\ge k_0$ and converges to $\xb$; that is, $\Ec$ lies in the basin of attraction of $\xb$.
\end{theorem}

\begin{proof}
That $O_c$ is open and has a unique component containing $\xb$ is shown in~\cite{RR2018}.

\emph{Forward invariance.} Let $x_k\in\Ec$ and $d_k=P_k\nabla f(x_k)$; we may assume $d_k\ne0$. Let $T=\{t\in[0,1]:x_k+sd_k\in\Ec\text{ for all }s\in[0,t]\}$ and $t^*=\sup T$. For $s<t^*$ the segment $[x_k,x_k+sd_k]$ lies in $\Ec$, where $\nabla f$ is $L$-Lipschitz, so~\ref{c:segment}, applied to the matrix $sP_k$ (whose step is this segment, and $\lambda_n(sP_k)\le\Lambda$), gives $f(x_k+sd_k)\ge f(x_k)>c$. Letting $s\rightarrow t^{*}$ from below, the point $y=x_k+t^*d_k$ lies in $\Ecb\subset D$ and $f(y)\ge f(x_k)>c$, so $y\notin\partial\Ec$ by H5; hence $y\in\Ec$ and $t^*\in T$. If $t^*<1$, then since $\Ec$ is open, a ball about $y$ lies in $\Ec$, so $t^*+\varepsilon\in T$ for some $\varepsilon>0$, contradicting the choice of $t^*$. Thus $t^*=1$, that is, $[x_k,x_{k+1}]\subset\Ec$. By induction $x_k\in\Ec\subset D$ for all $k\ge k_0$, and every step lies where the Lipschitz bound holds.

\emph{Convergence.} If some $x_k=\xb$ we are done. Otherwise $V_k$ is nonincreasing and nonnegative, with limit $\alpha\ge0$. If $\alpha>0$, choose $\delta>0$ with $V<\alpha$ on $B_\delta(\xb)$; every iterate lies in the compact set $\Ecb\setminus B_\delta(\xb)$, where $|\nabla f|^2$ has a positive minimum by H6, and~\ref{c:strict} with H3$_\Delta$ forces $V_k\to-\infty$, a contradiction. So $V_k\to0$, and $x_k\to\xb$ as in the proof of Theorem~\ref{thm:stab}(ii).
\end{proof}

In continuous time the segment argument is unnecessary, because a trajectory cannot leave $\Ec$ without crossing $\partial\Ec$. Here the bound $L\Lambda\le2-\eta$ in S3 ensures that $f$ does not dip anywhere along a step, so a step cannot jump the valley separating $\Ec$ from another component of $\{f>c\}$.

\begin{figure}[!ht]
\centering
\includegraphics[width=0.6\textwidth]{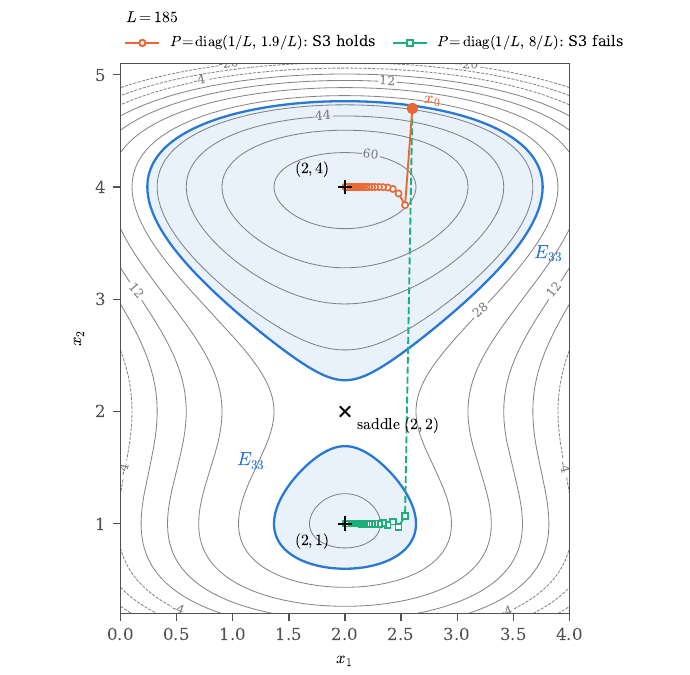}
\caption{Level curves of $f$ in Example~\ref{ex:basin}, with the two components of $E_{33}$ shaded and outlined; $+$ marks the peaks and $\times$ the saddle. Both trajectories start at $x_0=(2.6,4.70)$ and use $L=185$. Circles: $P_k=\operatorname{diag}(1/L,\,1.9/L)$, which satisfies S3; the first step overshoots the peak in $x_2$, but every iterate stays in $E_{33}$ and the iterates converge to $(2,4)$. Squares: $P_k=\operatorname{diag}(1/L,\,8/L)$, which violates S3; the first step (dashed) jumps the valley into $E_{33}$ about $(2,1)$, and the iterates instead converge to $(2,1)$ not $(2,4)$.}
\label{fig:basin}
\end{figure}

\begin{example}\label{ex:basin}
For Example~3.1 of~\cite{RR2018}, $f(x_1,x_2)=96x_2-84x_2^2+28x_2^3-3x_2^4-10(x_1-2)^2$ has peaks at $(2,1)$ and $(2,4)$ and a saddle at $(2,2)$. The sets $E_{33}$ about each peak lie in the corresponding basin of~\eqref{eq:Delta} once $\Lambda$ is below $2/L$, where $L$ is a Lipschitz constant for $\nabla f$ on $\bar E_{33}$. Because $\nabla^2f=\operatorname{diag}\bigl(-20,\,-168+168x_2-36x_2^2\bigr)$ depends only on $x_2$, $L$ is the maximum of $|\partial^2f/\partial x_2^2|$ over the $x_2$-range of $\bar E_{33}$.  Looking at the region about $(2,4)$ means that  $x_2 \in [2.279,4.764]$ (approximately) and $L\approx184.8$, attained at the top of the set, so we let $L=185$. In comparison, looking at the region about $(2,1)$ means that  $x_2 \in [0.599,1.691]$ (approximately) and $L\approx80.3$. With $P_k=\operatorname{diag}(1/L,\,1.9/L)\approx\operatorname{diag}(0.0054,\,0.0103)$, S3 holds with $\eta=0.1$, and an iterate started near the top of $E_{33}$ about $(2,4)$ overshoots the peak on its first step but stays in $E_{33}$ and converges to $(2,4)$ (Figure~\ref{fig:basin}). With $P_k=\operatorname{diag}(1/L,\,8/L)\approx\operatorname{diag}(0.0054,\,0.043)$, so that $L\Lambda=8$, the first step from the same point crosses the valley at the saddle and lands in $E_{33}$ about $(2,1)$, and the iterates subsequently converge to $(2,1)$. Note that the bound $\Lambda<2/L$ is conservative here in two ways. First, after its first step the second trajectory lies in $E_{33}$ about $(2,1)$, where $L\approx80.3$ and S3 would require $\Lambda<2/L\approx0.025$. Even so, with $\Lambda\approx0.043$ every later iterate stays in that set and the iterates converge to $(2,1)$, because linear stability at $(2,1)$, where $\partial^2f/\partial x_2^2=-36$, requires only $p_2<2/36\approx0.056$. Second, for constant diagonal $P=\operatorname{diag}(p_1,p_2)$, a step from $E_{33}$ about $(2,4)$ into $E_{33}$ about $(2,1)$ must decrease $x_2$ by at least $4.764-1.691\approx3.07$, while $|\partial f/\partial x_2|\le95.5$ on $\bar E_{33}$ about $(2,4)$; therefore $p_2\ge0.032$, about three times $2/L\approx0.011$. Because $\partial^2f/\partial x_2^2(2,4)=-72$ and $p_2>2/72\approx0.028$, such a $P$ already makes $(2,4)$ linearly unstable. The discussion in~\cite{RR2018} of the sets $E_{20}$, which violate H5 and H6, carries over unchanged.
\end{example}

\section{Necessity of the hypotheses}\label{sec:necessity}

Each row of Table \ref{tab:hypotheses} drops one hypothesis of Theorem~\ref{thm:stab} and keeps the others (because S2 implies S1, the first row drops both). The second row of the table deserves comment: As noted therein, Example~\ref{ex:quad} with $pa=2$ (Figure~\ref{fig:quad}b) keeps S1 but drops the margin $\eta$ in S2, and its trajectories are not asymptotically stable. Assuming S2 and H3$_\Delta$ leads to asymptotic stability (Theorem \ref{thm:stab}(ii)). However, Corollary \ref{cor:weighted} proves a weaker hypothesis than the combination of S2 and H3$_\Delta$ produces uniform asymptotic stability (Example \ref{ex:nomargin}). Therefore, S2 is sufficient but not necessary.

\begin{table}[!ht]
\tbl{Effects of dropping various hypotheses on the dynamics of~\eqref{eq:Delta}.\label{tab:hypotheses}}
{\begin{tabular}{@{}lll@{}}
\toprule
Hypothesis dropped & Example & What fails\\
\midrule
S1 ($L\Lambda\le2$) & \ref{ex:quad}, $pa>2$ & not stable; $|x_k|\to\infty$\\
margin $\eta$ in S2 & \ref{ex:quad}, $pa=2$ & stable, period-2 orbit, not asymptotically stable\\
H3$_\Delta$ & \ref{ex:EC} & stable, but $x_1$ converges to the wrong value\\
H2 & \ref{ex:rings} & stable, iterates trapped between rings of equilibria\\
\bottomrule
\end{tabular}}
\end{table}

\begin{figure}[!ht]
\centering
\includegraphics[width=0.6\textwidth]{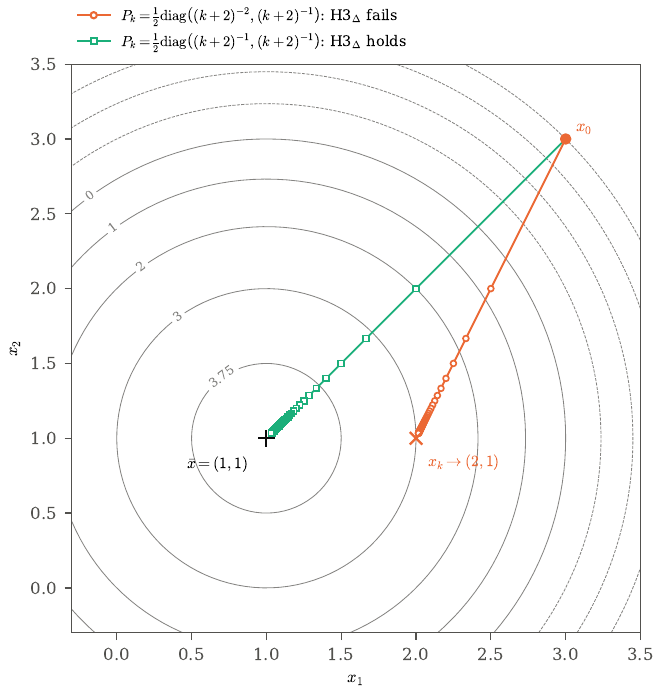}
\caption{Level curves of $f$ in Example~\ref{ex:EC}; $+$ marks the maximum $\bar x=(1,1)$. Both trajectories start at $x_0=(3,3)$ and have $\Lambda=1/4$. Circles: $P_k=\tfrac12\operatorname{diag}\bigl((k+2)^{-2},(k+2)^{-1}\bigr)$, for which H3$_\Delta$ fails; $x_{2,k}\to1$ but $x_{1,k}\to2$, so the iterates converge to $(2,1)$ ($\times$) rather than $\bar x$. Squares: $P_k=\tfrac12\operatorname{diag}\bigl((k+2)^{-1},(k+2)^{-1}\bigr)$, for which H3$_\Delta$ holds; $e_k=e_0/(k+1)\to0$ and the iterates converge to $\bar x$.}
\label{fig:EC}
\end{figure}

\begin{example}[The eigenvalue condition is needed]\label{ex:EC}
This is the discrete form of Example~2.1 of~\cite{RR2018}. Let $f(x_1,x_2)=4-(x_1-1)^2-(x_2-1)^2$, so $L=2$ and $\bar x = (1,1)$. Next let $P_k=\tfrac12\operatorname{diag}\bigl((k+2)^{-2},(k+2)^{-1}\bigr)$. Then the largest eigenvalue of $P_k$ occurs when $k=0$ giving $\Lambda=1/4$ and S2 holds with $\eta=3/2$. If we write out ($\Delta$) for this modified gradient system and substitute $e_k=x_k-(1,1)$, we get:

\begin{align*}
  e_{1,k+1}&=\bigl(1-(k+2)^{-2}\bigr)e_{1,k} \implies e_{1,k} = e_{1,0} \prod_{j=0}^{k-1}\bigl(1-(j+2)^{-2}\bigr)=e_{1,0}\,\frac{k+2}{2(k+1)},\\ e_{2,k+1}&=\bigl(1-(k+2)^{-1}\bigr)e_{2,k} \implies e_{2,k} = e_{2,0}\prod_{j=0}^{k-1}\bigl(1-(j+2)^{-1}\bigr)=\frac{e_{2,0}}{k+1}.
\end{align*}

\noindent Starting at $k_0=0$ gives $e_{1,k}\to e_{1,0}/2$ while $e_{2,k}\to0$, so $x_k\to(1+e_{1,0}/2,\,1)\ne\bar x$ whenever $e_{1,0}\ne0$ (Figure \ref{fig:EC}). Note that H3$_\Delta$ has been dropped because $\sum_k\lambda_1(P_k) = \tfrac12\sum_k (k+2)^{-2} = (\pi^2 - 6)/12<\infty$.  In Example~2.1 of~\cite{RR2018}, $P(t)=\operatorname{diag}\bigl((t+1)^{-2},(t+1)^{-1}\bigr)$; here we add the factor $\tfrac12$ and shift to $k+2$. The shift to $k+2$ avoids a zero factor at $k=0$, which would send $x_1$ to $\bar x$ in a single step.
\end{example}

\begin{figure}[!ht]
\centering
\includegraphics[width=\textwidth]{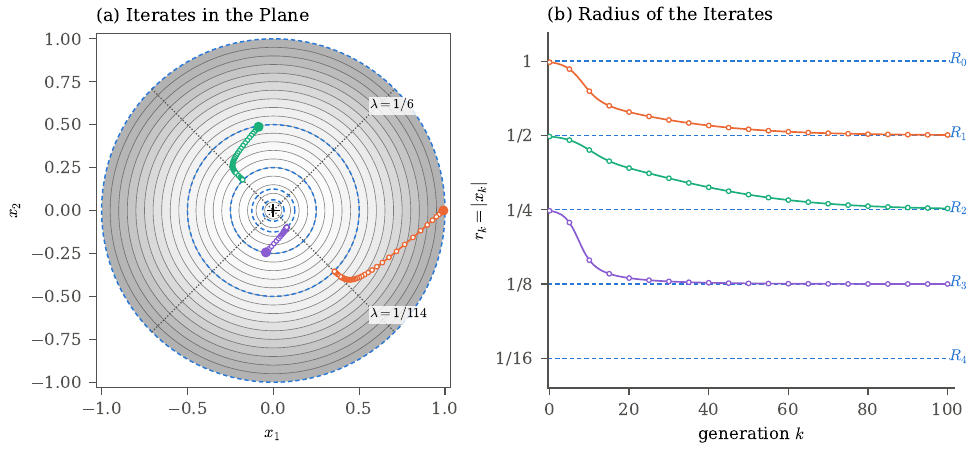}
\caption{Iterates of Example~\ref{ex:rings}. (a) Level curves of $f(x)=p(|x|)$ on the unit disk; dashed circles are the rings of critical points $r=R_n=2^{-n}$, $n=0,\dots,4$. $+$ marks the maximum $\bar x=0$. In this example $P_k=\frac1{114}\bigl(\begin{smallmatrix}10&9\\9&10\end{smallmatrix}\bigr)$, and dotted lines are the eigendirections of $P$ (eigenvalues $1/6$ and $1/114$). Each trajectory starts just inside a ring $R_n$ ($r_0=0.99R_n$, $n=0,1,2$). The component along the fast direction decays within a few steps, and the iterates then move along the slow direction, but each trajectory stops at the next ring $R_{n+1}$. (b) The radius $r_k=|x_k|$ on a logarithmic scale: $r_k$ decreases to $R_{n+1}$ and never crosses it, so the iterates do not converge to $\bar x$.}
\label{fig:rings}
\end{figure}

\begin{example}[An isolated critical point is needed]\label{ex:rings}
This is the discrete form of Example~2.2 of~\cite{RR2018}. Let $p$ be the even $C^1$ spline on $[-1,1]$ constructed there, with critical points at the radii $R_n=2^{-n}$, and let $f(x)=p(|x|)$ on the closed unit disk. Between consecutive critical radii $|p''|\le6$, and $|p'(r)/r|\le3/2$, so $\nabla f$ is Lipschitz with $L=6$ even though $f$ is not $C^2$. Take
\[
  P_k=P=U\operatorname{diag}\bigl(\tfrac16,\tfrac1{114}\bigr)U^{\mathsf T}=\frac1{114}\begin{pmatrix}10&9\\9&10\end{pmatrix},
\]
where $U$ is the rotation by $\pi/4$. Then $\Lambda=1/6$, so $L\Lambda=1$ and S2 holds with $\eta=1$, and H3$_\Delta$ holds because $\lambda_1(P)=1/114>0$. Because $\nabla f(x)=\bigl(p'(r)/r\bigr)x$, the iterates satisfy $x_{k+1}=(I+c_kP)x_k$ with $c_k=p'(r_k)/r_k\in[-3/2,0]$. The eigenvalues of $I+c_kP$ lie in $[1+c_k\Lambda,1]\subset(0,1]$, so $r_k+\Lambda p'(r_k)\le r_{k+1}\le r_k$. If $R_{n+1}<r_k<R_n$, then $p'(R_{n+1})=0$ and the Lipschitz bound give $|p'(r_k)|\le6(r_k-R_{n+1})$, therefore
\[
  R_{n+1}\le r_k-6\Lambda(r_k-R_{n+1})\le r_{k+1}\le r_k<R_n .
\]
In other words, an iterate never crosses the ring $r=R_{n+1}$, so $\xb=0$ is stable but not asymptotically stable (Figure~\ref{fig:rings}). 

\end{example}

As noted in~\cite{RR2018}, H2 is missing from some textbook statements for gradient systems~\cite[p.~205]{HSD2004}. Additional structure can partly replace it: for real-analytic $f$ the \L{}ojasiewicz inequality forces the iterates of descent methods to converge to a single critical point~\cite{Absil2005}, though not necessarily to $\xb$. For Theorem~\ref{thm:basin}, Example~3.1 of~\cite{RR2018} already shows that H5 and H6 are needed, and Example~\ref{ex:quad} shows that a step with $L\Lambda>2$ can leave $\Ec$.

\FloatBarrier
\section{Discussion}\label{sec:discussion}

Previously, it has been show that for the continuous-time multivariate breeder's equation that an isolated local maximum of  the population fitness surface $f = \ln \bar W$, whose gradient function describes the forces of natural selection, is uniformly stable whenever $P(t)$ is positive semidefinite \cite{RR2018}. In the case of evolutionary systems, $P(t)$ is the time-varying ``$G$-matrix'' which is \emph{always} positive semidefinite. Here we show the same cannot be said for the discrete-time analog. However it is possible to recover similar results for such systems if we additionally impose a step-size bound. Specifically, if hypothesis S1 holds ($L\Lambda \le 2$), we recover uniform stability; furthermore, if hypothesis S2 holds ($L\Lambda \le 2 - \eta, \, \eta \in (0,2]$) together with a hypothesis (H3$_\Delta$) similar to the eigenvalue condition of \cite{RR2018} then  uniform asymptotic stability is recovered. In these hypotheses, $L$ is the Lipschitz constant of $\nabla f$, and $\Lambda$ bounds the largest eigenvalue of every $P_k$.

It should be emphasized that the Lipschitz constant is local. If $f\in C^2$, then $L$ can be taken arbitrarily close to $\|\nabla^2f(\xb)\|$ by shrinking $r_0$, so $\sup_k\lambda_n(P_k)<2/\|\nabla^2f(\xb)\|$ suffices for Theorem~\ref{thm:stab}. As an extreme example of this, if $f= -x^4$, then $\nabla^2f(\xb)=0$ and any bounded sequence $P_k$ gives stability. Of more interest is how $L$ grows with $\operatorname{dist}(x,\xb)$, which determines how large $\lambda_n(P_k)$ can be while stability is maintained. As discussed below, this has implications for the importance of variance-generating genetic processes (e.g., mutation, recombination) and for evolutionary algorithms used in optimization.

\subsection{Biological interpretation}
This work has been inspired by the study of evolutionary dynamics and thus it is natural to interpret our results biologically. Recall that $P_k=G_k$ is the genetic covariance matrix, and $\nabla f=\ln\bar W$ is the natural selection gradient \cite{Lande1976, LandeArnold1983}. Conditions S1 and S2 say that the additive genetic variance in every direction must be small relative to the curvature of log mean fitness near the peak. Within this regime the conclusions of~\cite{RR2018} remain: the evolution of $G$ may change the path and the speed of adaptation but cannot destabilize a fitness peak. Outside it (i.e.,when S1 fails), the mean phenotype can overshoot the optimum and oscillate or diverge (Example~\ref{ex:quad}). This is not the only mechanism by which an evolutionary optimum might not be achieved. For a different mechanism by which a fitness maximum can fail to attract see~\cite{Abrams1993}; it and similar papers study frequency-dependent selection where the fitness function $f$ changes in time.

Our findings with respect to the size of eigenvalues that produce stable dynamics relate to the concept of ``overshoot'' in Fisher's geometric model \cite{Fisher1930}. Fisher made a similar argument: a phenotype at distance $d$ from an optimum in $n$ dimensions is perturbed by a mutation of random direction and size $r$. Fisher argued that mutations of very small effect are beneficial about half the time, whereas mutations of large effect are almost always harmful, because they carry the phenotype past (overshoot) the optimum. Even a mutation pointing directly at the optimum improves fitness only if $r<2d$; a longer step lands on the far side, farther from the optimum than it started. Kimura~\cite{Kimura1983} noted that large beneficial mutations are more likely to fix, so mutations of intermediate size contribute most to adaptation, and Orr~\cite{Orr1998,Orr2005} showed that the factors fixed during an adaptive walk have an approximately exponential distribution of effects. These studies show that having mutations that produce small changes in phenotype (so-called genes of minor effect~\cite{Mather1941,Lande1983,Rockman2012}) is critical to achieving evolutionary optima.

Example~\ref{ex:quad} directly relates to Fisher's geometric model. The response $P_k\nabla f(x_k)$ always points uphill, so, unlike a random mutation, a short enough response always improves $f$ (Lemma~\ref{lem:ascent}). Its length, however, is governed by $\lambda_n(P_k)$, the genetic variance in the direction of greatest variance. In Example~\ref{ex:quad} the step from $x_k$ has length $pa|x_k|$: it crosses the optimum when $pa>1$, and it lands farther away than it started, exactly Fisher's condition $r>2d$, when $pa>2$ (Figure~\ref{fig:quad}). Large genetic variance thus plays the role of a mutation of large effect, and conditions S1 and S2 are the requirement that the per-generation response never be ``large'' in Fisher's sense relative to the curvature of $\ln\bar W$. Theorems~\ref{thm:stab} and~\ref{thm:basin} show that, under this requirement (together with the eigenvalue condition for convergence, H3$_\Delta$), the conclusions of~\cite{RR2018} hold however $G$ evolves.

\subsection{Conservative isotropic bound} 
The classical model of Gaussian stabilizing selection shows that the S1 bound is conservative for some fitness models. Suppose individual fitness is $W(z)=\exp\bigl(-\tfrac12(z-\theta)^\T\Omega^{-1}(z-\theta)\bigr)$ and phenotypes are normally distributed with covariance $P_z=G+E$, where $E$ is the environmental covariance. Then, up to an additive constant, $\ln\bar W(\bar z)=-\tfrac12(\bar z-\theta)^\T(\Omega+P_z)^{-1}(\bar z-\theta)$~\cite{Lande1979}, so~\eqref{eq:lande} with constant $G$ and $P_z$ is the linear recursion
\[
  \bar z_{k+1}-\theta=\bigl(I-G(\Omega+P_z)^{-1}\bigr)(\bar z_k-\theta).
\]
The matrix $G(\Omega+P_z)^{-1}$ is similar to $(\Omega+P_z)^{-1/2}G(\Omega+P_z)^{-1/2}$, and $G\preceq P_z\prec\Omega+P_z$, so its eigenvalues lie in $[0,1)$ and the optimum is always approached without overshoot. In one dimension this is the familiar fact that the response per generation is the heritable fraction $G/(\omega^2+P_z)<h^2\le1$ of the distance to the optimum (with $\Omega = \omega^2$). In several dimensions, however, S1 can fail for this model: $L=1/\lambda_1(\Omega+P_z)$, and $\lambda_n(G)$ may exceed $2\lambda_1(\Omega+P_z)$ when selection is strong in one direction while genetic variance is large in another. Overshoot therefore requires log-fitness surfaces that are more sharply curved than the Gaussian model allows relative to $G$, and \emph{the isotropic bound S1 should be read as a sufficient condition that is tight only in the worst direction.}  

\subsection{Evolutionary optimization algorithms}
Evolutionary algorithms are now widely used for optimization in machine learning and elsewhere~\cite{HansenOstermeier2001,Wierstra2014,Salimans2017}, and the breeder's equation has been used explicitly to predict their progress~\cite{MuhlenbeinSV1993}. Our results apply directly to the mean of a population with a Gaussian search distribution. Suppose candidate solutions are drawn as $x\sim N(m_k,\Sigma_k)$, each is assigned a positive fitness $W(x)$, and the next mean is the fitness-weighted mean of the sample. In the limit of a large sample, Stein's identity gives
\[
  m_{k+1}-m_k=\frac{\operatorname{Cov}(x,W)}{\bar W_k(m_k)}=\Sigma_k\,\nabla\ln\bar W_k(m_k),\qquad \bar W_k(m)=\int W(x)\,\phi_{\Sigma_k}(x-m)\,dx,
\]
which is the same argument as Lande’s derivation for Gaussian phenotypes~\cite{Lande1979} and is system~\eqref{eq:Delta} with $P_k=\Sigma_k$ and $f_k=\ln\bar W_k$. With a learning rate $\alpha$ and a search distribution $N(m_k,\sigma_k^2C_k)$, the same calculation gives $P_k=\alpha\sigma_k^2C_k$; this is the Gaussian-smoothed gradient step used by evolution strategies in reinforcement learning~\cite{Salimans2017,NesterovSpokoiny2017}. Natural evolution strategies and CMA-ES (covariance matrix adaptation evolution strategy) update the mean by natural-gradient steps of the same form, applied to a rank-based transformation of the objective~\cite{Wierstra2014,Ollivier2017}.

In this setting the ``genetic variance'' is the search covariance, and it is shaped by every design choice: the spread of the initial population, the mutation distribution, recombination, the selection scheme, and the loss of diversity through sampling (an algorithm's analog of genetic drift). The hypotheses of Theorems~\ref{thm:stab} and~\ref{thm:basin} translate into guidance for these choices in three ways:
\begin{itemize}
  \item \emph{Upper bound.} If $\lambda_n(P_k)$ exceeds $2/L$, the mean can overshoot a peak, oscillate about it or diverge (Example~\ref{ex:quad}), however many generations are run. A widely initialized population or a large mutation variance is therefore not harmless near an optimum. That an overly large mutation strength makes progress negative is well known for evolution strategies~\cite{Beyer2001}, and step-size and covariance adaptation~\cite{HansenOstermeier2001} can be read as mechanisms that keep $P_k$ within the bounds S1--S3.
  \item \emph{Lower bound.} If the variance collapses so quickly in some direction that $\sum_k\lambda_1(P_k)<\infty$, the mean can stall short of the optimum in that direction (Example~\ref{ex:EC}). This is a precise form of premature convergence. Mutation that keeps $\lambda_1(P_k)$ bounded away from zero, or decays it no faster than $1/k$, satisfies H3$_\Delta$, in the same way that the classical step-size conditions of stochastic approximation require $\sum_k\alpha_k=\infty$~\cite{RobbinsMonro1951}.
  \item \emph{Basins.} Under S3, a population whose mean starts in $\Ec$ cannot leave it (Theorem~\ref{thm:basin}), so a small search covariance commits the algorithm to the nearest peak. Conversely, a covariance that violates S3 can allow the mean to jump the valley between peaks (Figure~\ref{fig:basin}). Our bounds thus quantify one side of the familiar trade-off between exploitation and exploration.
\end{itemize}
Two cautions apply. First, $\bar W_k$ depends on $\Sigma_k$, so when the search covariance changes the smoothed objective changes as well, whereas our theorems assume a fixed $f$; the results apply directly when $\Sigma_k$ is constant or when $\bar W_k$ changes slowly. Second, a finite sample adds noise to each update, which leads to the stochastic version of~\eqref{eq:Delta} raised in the open problems below. It is also worth noting that rank-based selection and discrete encodings are outside the present framework. Rank-based selection may nonetheless be brought within reach through quantitative-genetic models of threshold traits, in which a discrete or ordinal phenotype is expressed when an underlying, normally distributed breeding value exceeds one or more thresholds~\cite{Wright1934,DempsterLerner1950,Falconer1965,Roff1996}; selection on the observed categories then acts on the mean breeding value through an equation of the form~\eqref{eq:Delta}, so these models offer a route for extending the present results to ranked evolutionary algorithms.

\subsection{Conclusions and open problems}

Our work gives clear bounds as to when evolutionary optima should be attained in systems where discrete generations occur and the traits under investigation are continuous. These findings directly apply to the study of $G$-matrices, their evolution, and the resulting phenotypic dynamics. Furthermore, this research adds insight into ways to improve the implementation of genetic algorithms; this is particularly notable given the ever-increasing use of genetic algorithms to accelerate optimization routines~\cite{Salimans2017,HansenOstermeier2001}. 

Our work also poses several potential open problems to be addressed by future researchers. First of these is the search for an anisotropic bound. Because $P_k$ varies in time, the bound $L\Lambda \le 2$ is conservative and a more lenient condition might be found via a joint spectral radius problem \cite{Jungers2009}. Second, Example \ref{ex:rings}, where $\eta = 1$, shows trapping of trajectories between rings of critical points; it is unclear when such rings will trap the trajectory and when the trajectory might proceed to some other ring or eventually to $\xb$. Third, the addition of genetic drift (stochasticity) may have important implications: it could either lead to some concentration near $\xb$ or produce some expected time to leave the basin of attraction. (This is also important to the understanding of genetic algorithms.) Lastly, it is natural to ask how small the step-size in \eqref{eq:Delta} must be to converge to the basin of attraction results of \cite{RR2018}.

\section*{Disclosure statement}
The author reports there are no competing interests to declare.

\section*{Data availability statement}
No new data were generated in this study. Machine-checked proofs, in Lean~4 with the Mathlib library, of Lemma~\ref{lem:ascent}, Theorem~\ref{thm:stab}, Corollary~\ref{cor:weighted} and Theorem~\ref{thm:basin} are openly available on GitHub at \url{https://github.com/magbenji/discrete-modified-gradient-lean} and archived on Zenodo at \url{https://doi.org/10.5281/zenodo.24118844} (version~1.0).

\section*{Funding}

None.

\bibliographystyle{tfs}
\bibliography{references}

\end{document}